\documentclass{elsarticle}
\usepackage[latin1]{inputenc}
\usepackage{amscd}
\usepackage{amsmath,amssymb,amsfonts}
\usepackage{syntonly}

\usepackage{graphics, graphicx}
\usepackage{color}
\usepackage{epsfig}
\usepackage{mathrsfs, amsfonts}

\newcommand{\R}[0]{\mathbb{R}}
\newcommand{\C}[0]{\mathbb{C}}

\renewcommand{\S}{\mathcal{S}}

\newtheorem{thm}{Theorem}[section]
\newtheorem{prop}{Proposition}

\newtheorem{lm}[thm]{Lemma}
\newtheorem{crl}[thm]{Corollary}
\newtheorem{claim}{Claim}[thm]
\newdefinition{nt}[thm]{Notation }
\newdefinition{qu}[thm]{Question }
\newdefinition{rem}[thm]{Remark }
\newdefinition{tdef}[thm]{Theorem-Definition }
\newdefinition{df}[thm]{Definition }
\newdefinition{ex}[thm]{Example }

\newproof{pf}{Proof}
\newproof{pf1}{Proof of claim \ref{banale}}
\newproof{pf2}{Proof of claim \ref{claim2}}
\newproof{pft}{Proof of Theorem \ref{crescej}}
\newproof{pfsk}{Sketch of proof}
\newproof{pfcl}{Proof of the claim}

\renewcommand{\AA}{\mathcal{A}}
\newcommand{\PP}{\mathcal{P}}
\newcommand{\FF}{\mathcal{F}}

\newcommand{\LL}{\mathcal{L}}
\newcommand{\MM}{\mathcal{M}}
\renewcommand{\AA}{\mathcal{A}}

\newcommand{\CC}{\mathcal{C}}

\newcommand{\pol}{\lhd} 
\newcommand{\fl}{\prec} 
\newcommand{\fleq}{\preceq} 
\newcommand{\fgeq}{\succeq} 
\newcommand{\fg}{\succ} 
\newcommand{\Ha}{\widetilde{H}}

\newcommand{\UU}{\mathcal{U}}
\newcommand{\La}{\widetilde{L}}
\newcommand{\lkj}{\La^{k,j}}

\newcommand{\op}{\textup{op}}

\newcommand{\codim}{\textup{codim}}

\newcommand{\XX}{\mathcal{X}}
\newcommand{\YY}{\mathcal{Y}}
\newcommand{\HH}{\mathcal{H}}

\usepackage{ifpdf}
\ifpdf
\usepackage[%
  pdftitle={Instructions for use of the document class
    elsart},%
  pdfauthor={Simon Pepping},%
  pdfsubject={The preprint document class elsart},%
  pdfkeywords={instructions for use, elsart, document class},%
  pdfstartview=FitH,%
  bookmarks=true,%
  bookmarksopen=true,%
  breaklinks=true,%
  colorlinks=true,%
  linkcolor=blue,anchorcolor=blue,%
  citecolor=blue,filecolor=blue,%
  menucolor=blue,pagecolor=blue,%
  urlcolor=blue]{hyperref}
\else
\usepackage[%
  breaklinks=true,%
  colorlinks=true,%
  linkcolor=blue,anchorcolor=blue,%
  citecolor=blue,filecolor=blue,%
  menucolor=blue,pagecolor=blue,%
  urlcolor=blue]{hyperref}
\fi

\makeatletter
\def\elsartstyle{%
    \def\normalsize{\@setfontsize\normalsize\@xiipt{14.5}}
    \def\small{\@setfontsize\small\@xipt{13.6}}
    \let\footnotesize=\small
    \def\large{\@setfontsize\large\@xivpt{18}}
    \def\Large{\@setfontsize\Large\@xviipt{22}}
    \skip\@mpfootins = 18\p@ \@plus 2\p@
    \normalsize
}
\@ifundefined{square}{}{}
\makeatother

\begin{document}

\begin{frontmatter}

\title{Combinatorial polar orderings\\ and recursively orderable arrangements}

\author{Emanuele Delucchi\fnref{msri}}
\fntext[msri]{The first author was partially supported by a postdoctoral fellowship of the Mathematical Sciences Research Institute at Berkeley, CA, USA.}
\address{Department of Mathematics, State  University of New York, PO Box 6000,
Binghamton, NY 13902-6000, USA. }
\ead{delucchi@binghamton.edu}

\author{Simona Settepanella}
\address{Laboratory of Economics and Management,
Scuola Superiore Sant'Anna,
Piazza Martiri della Libert\`a 33,
56127 Pisa, Italy.} 
\ead{s.settepanella@sssup.it}


\begin{abstract} Polar orderings arose in recent work of Salvetti and the second
  author on minimal CW-complexes for complexified hyperplane
  arrangements. We study the combinatorics of these orderings in the classical framework of oriented matroids, and reach thereby
  a weakening of the conditions required to actually determine such
  orderings. A class of arrangements for which the construction of
  the minimal complex is particularly easy, called {\em recursively orderable} arrangements, can therefore be
  combinatorially defined. We initiate the study of this class, giving
  a complete characterization in dimension 2 and proving that every
  supersolvable complexified arrangement is  recursively orderable.
\end{abstract}
\end{frontmatter}
\section*{Introduction}

One of the main topics in the theory of arrangements of hyperplanes is
the study of the topology of the complement of a set of hyperplanes in
complex space. The special case of {\em complexified arrangements}, where the hyperplanes have real
defining equations, is very interesting in its own as it allows
a particularly explicit combinatorial treatment. Indeed, when dealing
with complexified arrangements one can rely
on the {\em Salvetti complex}, a regular CW-complex that can be
constructed entirely in terms of the oriented matroid of the real arrangements and is a deformation retract
of the complement of the complexified arrangement \cite{salvetti}.

A general fact about complex arrangement's complements is that they
are {\em minimal spaces} (i.e., they carry the homotopy type of a
CW-complex where the number of cells of any given dimension equals the
rank of the corresponding homology group), as was proved by Dimca and
Papadima \cite{dimcapapa} and, independently, by Randell \cite{randell} using Morse theoretical arguments. Again, in
the complexified case the topic allows an explicit treatment:
as shown in \cite{salvettisette,delu}, one can exploit discrete Morse theory on the
Salvetti complex to construct a discrete Morse vector field that
allows to collapse every `superfluous' cell and thus produces an
explicit instance of the minimal complex whose existence was predicted
in \cite{dimcapapa,randell}.

The approach taken by Salvetti and the second author in
\cite{salvettisette} to construct the discrete Morse vector field
relies on the choice of a so-called {\em generic flag} and on the
associated {\em polar ordering} of the faces of the real
arrangement. Once this polar ordering is determined, the description
of the vector field and of the obtained minimal complex is quite
handy, e.g. yielding an explicit formula for the algebraic boundary
maps.

But the issue about actually {\em constructing} such a polar ordering
for a given arrangement remains. This motivates the first part of our
work, where we give a fully combinatorial characterization of a whole class of
total orderings of the faces of a complexified arrangement that can be
used as well to carry out the construction of the very same discrete vector
field described in \cite{salvettisette}. Our {\em combinatorial polar
  orderings} still require a flag of general position subspaces as a starting
point, but does not need this flag to satisfy the requirements that
are requested from a {\em generic flag} in the sense of
\cite{salvettisette}. Our construction builds upon the
concept of {\em flipping} in oriented matroids, letting a
pseudohyperplane `sweep' through the arrangement instead of `rotating'
it around a fixed codimension $2$ subspace as in \cite{salvettisette}
(see our opening section for a review of the concepts).

Once the (combinatorial) polar ordering is constructed, one has to
figure out the discrete vector field and follow its gradient paths to
actually construct the minimal complex. Although the `recipe' is
fairly straightforward, this task soon becomes very
challenging. For instance, this was accomplished in \cite{salvettisette} for the family of
 real reflection arrangements of Coxeter type $A_n$. The key fact
 allowing one to carry out the construction in these cases is that the
 general flag can be set so that the associated polar orderings enjoy
 a special technical property (see Definition \ref{df:maindef}) that
 keeps the complexity of computations down to a reasonable level.

Thus it is natural to ask whether this property is shared by other
arrangements. Since the obtained discrete vector fields are the same, it turns out that instead of restricting to `actual' polar
orderings, it is natural to work in our broader combinatorial setting,
and say that an arrangement is {\em recursively orderable} if it admits a
combinatorial polar ordering that satisfies the same technical
property that made computations feasible for the $A_n$ arrangements.

In the second part of this work we initiate the study of recursively orderable
arrangements. We reach a complete characterization
of this property for arrangements of lines. Trying to generalize
the property to the three major classes of arrangements to which $A_n$
belongs, we prove that every supersolvable arrangement is
recursively orderable. Indeed, the required recursive ordering can be recovered basically from the standard decomposition into ``blocks'' (i.e., modular flats) of supersolvable arrangements. On the other hand, not every reflection arrangement is recursively orderable. As what concerns asphericity, already in dimension $3$ there is a
recursively orderable arrangement that is not $K(\pi,1)$.
We believe that the class of recursively orderable
arrangements still bear some combinatorial and topological interest,
and deserve further study.\\

The paper starts with a section that gives some theoretical background
and reviews the different techniques needed later on.

Then the first part of the actual work is dedicated to the
combinatorial study of polar orderings. We begin by explaining the setup and
the required notation for handling with flippings of affine oriented
matroids. Then, in Section \ref{special} we give some
characterization of the valid sequences of flippings that allow a
pseudohyperplane to sweep across an affine arrangement, and call these
{\em special orderings} of the points of the arrangement. A key fact
in this section is how special orderings of the
points of the arrangement {\em induced on the moving pseudohyperplane}
behave after each ``move'' of the pseudohyperplane. In this view, the genericity condition on the general flag of
\cite{salvettisette} ensures that every step in the sequence of
flippings leads to a realizable oriented matroid, on which a polar ordering can be defined with the same geometric construction. Now, the contraction of the arrangement $\AA$ to our moving pseudohyperplane may not in general give rise to a realizable oriented matroid. 
However, we can prove that at each step in our construction the contractions that have to be performed lead to configurations that, although not realizable, admit a `sweeping' as above. This fact is proved using the theory of {\em oriented matroid programs} (see Definition \ref{teodef}). Indeed, an oriented matroid program is  an affine oriented matroid with a distinguished element, and it is called `Euclidean' if and only if the (pseudo-)hyperplane corresponding to the distinguished element can be `swept' through the whole affine oriented matroid. In our case (Remark \ref{OMP_expl}) we check an equivalent caracterization of this property established by Fukuda (see \cite[Chapter 10]{BLSWZ} for reference).

 In Section \ref{combinatorial} we then associate a {\em combinatorial polar ordering} to every set of one special ordering
for every one of the sections of the arrangement induced on a flag of
generic subspaces. To prove that this definition indeed makes sense,
Section \ref{polar} shows that every combinatorial polar ordering can
be obtained from a `genuine' polar ordering by a sequence of moves,
called {\em switches}, that do not affect the induced discrete vector
field. Thus every combinatorial polar ordering induces a discrete
Morse function with a minimum possible number of critical cells, and
leads to a minimal complex for the arrangement's complement
(Proposition \ref{aaa}).

The second part of the work, as said, is devoted to recursively orderable 
arrangements. The definition is given in Section \ref{the} along with
some basic facts. Section \ref{follow} studies the $2$-dimensional case,
leading, with Theorem \ref{thm_2d}, to a necessary and sufficient
condition for an arrangement of lines to be recursively orderable. We close this
paper with Section \ref{supersolvable}, where we prove that every
supersolvable arrangement is recursively orderable.\\

We are pleased to thank Mario Salvetti for his helpful and encouraging advice. We also are grateful to Laura Anderson for the opportunity to discuss with her some issues of a first version of this paper.

\numberwithin{thm}{subsection}

\addtocounter{section}{18}
\renewcommand{\thesection}{\Alph{section}}

\section*{Review}

\subsection{Topology and combinatorics of complexified arrangements.}
Let $\AA$ be an essential affine hyperplane
arrangement in $\R^d$, i.e., a set of affine real hyperplanes whose minimal nonempty intersections are points. Let $\FF$ denote the set of closed strata of the induced stratification of $\R^d$. It is customary to endow $\FF$ with a partial ordering $\fleq$ given by reverse inclusion of topological closures. The elements of $\FF$ are called {\em faces} of the arrangement. Their closures are polyhedral subsets of $\R^d$ and therefore we will adopt the corresponding terminology; given $F\in\FF$, the {\em faces of $F$} are the polyhedral faces of the closure of $F$, and consistently a {\em facet} of $F$ is any maximal face in its boundary. The poset $\FF$ is ranked by the {\em codimension} of the faces. The connected components of $\mathbb{R}^d\setminus \AA$, corresponding to  elements of $\FF$ of maximal dimension, are called {\em chambers}. 
For any $F\in\FF$ let $|F|$ denote the affine subspace spanned by $F$, called the {\em support} of $F$, and set 
$$\AA_F\ :=\ \{ H \in \AA\ :\ F\subset H\}.$$

Mario Salvetti \cite{salvetti} constructed a regular CW-complex 
$\S(\AA)$ (denoted just by $\S$ if no misunderstanding about the arrangement can arise) that is a deformation retract of
$$\MM(\AA):=\ \C^d \setminus \bigcup_{H\in A}\ H_{\C},$$
the complement of the complexification of $\AA$.

The  $k$-cells of $\S$ bijectively correspond to pairs
$[C\fleq F]$ where $\codim(F)=k$ and $C$ is a chamber. A cell $[C_1\fleq F_1]$ is in the boundary of $[C_2\fleq F_2]$ if $F_1 \fl F_2$ and the chambers $C_1$, $C_2$ are contained in the same chamber of $\AA_{F_2}$.
\subsection*{Discrete Morse theory.} A combinatorial version of Morse theory that is particularly well-suited for working on regular CW-complexes was formulated by Forman \cite{ForDM}. Here we outline the basics of Forman's construction, and we point to the book of Kozlov \cite{Kozlov} for a broader introduction and a more recent exposition of the combinatorics of this subject.\\

\newcommand{\KK}{\mathcal{K}}

\begin{df} Let $K$ be a locally finite regular $CW$-complex and $\KK$ denote the set 
of cells of $K$, ordered by inclusion. A \textit{discrete Morse function} on $K$ is a function
$f:\KK \longrightarrow \R$ such that
$$\begin{array}{cccccc} (i) && \sharp \{\tau^{(p+1)}>\sigma^{(p)}\ |\ f(\tau^{(p+1)})
\leq f(\sigma^{(p)})\} & \leq & 1& \\
(ii) &&  \sharp \{\tau^{(p-1)} < \sigma^{(p)} \ |\ f(\sigma^{(p)})
\leq f(\tau^{(p-1)}) \} & \leq & 1&\\
\end{array}$$
for all cells $\sigma^{(p)}\in\KK$ of dimension $p$.

 Moreover, $\sigma^{(p)}$ is a \textit{critical cell
of index $p$} if both sets are empty. Let  $m_p(f)$ denote the number of critical cells of $f$ of index
$p$. 
\end{df}

This setup is a discrete analogue of classical Morse theory in the following sense.

\begin{thm}[\cite{ForDM}, see also \cite{Kozlov}] If $f$ is a discrete Morse function on the regular $CW$-complex $K$, then $K$ is homotopy equivalent to a 
$CW$-complex with exactly $m_p(f)$ cells of dimension $p$.
\end{thm}

\begin{df} Let $f$ be a discrete Morse function on a $CW$-complex $K$. The discrete gradient vector field $V_f$ of $f$ is:
$$V_f\ =\ \{(\sigma^{(p)},\tau^{(p+1)}) | \sigma^{(p)}>\tau^{(p+1)},\,f(\tau^{(p+1)}) \leq f(\sigma^{(p)})\}.$$\end{df}

By definition of Morse function, each cell belongs to at most one
pair of $V_f$. So $V_f$ is a matching of the edges of the Hasse diagram of $\FF$ and the critical cells are precisely the non-matched elements of $\KK$. Because $f$ is a discrete Morse function, there cannot be any cycle in $\FF$ that alternates between matched and unmatched edges - such a matching is called {\em acyclic}. The following is a crucial combinatorial property of discrete Morse functions.

\begin{thm}[\cite{Kozlov}]\label{acgrad} For every acyclic matching $M$ of $\KK$ there is a discrete Morse function $f$ on $K$ so that $M=V_f$.  Thus, discrete Morse functions on $K$ correspond to  acyclic matchings of the Hasse diagram of $\KK$.
\end{thm}

\subsection{Polar ordering and polar gradient.} 

Salvetti and the second author introduced {\em polar orderings} of real hyperplane arrangements in \cite{salvettisette} as the basic tool for the construction of minimal models for $\MM(\AA)$. The construction starts by considering the polar coordinate system induced by any  \textit{generic flag} with respect to the given arrangement $\AA\subset\mathbb{R}^d$, i.e., a flag $\{V_i\}_{i=0,\ldots,d}$ of affine subspaces in general position, such that $\dim(V_i)=i$ for every $i=0,\ldots,d$ and such that `the polar coordinates $(\rho,\theta_1,\ldots,\theta_{d-1})$ of every point in a bounded face of $\AA$ satisfy $\rho>0$ and $0< \theta_i < \pi/2$, for every $i=1,\ldots,d$' (see \cite[Section 4.2]{salvettisette} for the precise description). The existence of such a generic flag is not trivial (\cite[Theorem 2]{salvettisette}). Every face $F$ is labeled by the coordinates of the point in its closure that has lexicographically least polar coordinates. 

The {\em polar ordering} associated to a generic flag is the total order $\pol$ on $\FF$ that is obtained by ordering the faces lexicographically according to their labels. This extends the order in which $V_{d-1}$ intersects the faces while rotating around $V_{d-2}$. If two faces share the same label - thus, the same minimal point $p$ -, the ordering is determined by the general flag induced on the copy of $V_{d-1}$ that is rotated `just past $p$' and the ordering it generates by induction on the dimension (see \cite[Definition 4.7]{salvettisette}).

The main purpose of the polar ordering is to define a discrete Morse function on the Salvetti complex, which, by Theorem \ref{acgrad}, amounts to specifying an acyclic matching $\Phi$ on the poset of cells of $\S$ that is called the {\em polar gradient}. The original definition of $\Phi$ is by induction in the dimension of the subspace $V_k$ containing the faces \cite[Definition 4.6]{salvettisette}. For the sake of brevity let us here define $\Phi$ through an equivalent description that is actually the one we will use later (compare Definition \ref{Fi})

\begin{df}[Compare Theorem 6 of \cite{salvettisette}] For any two faces $F_1,F_2$ with $F_1\fl F_2$, $\codim(F_1)=\codim(F_2)-1$ and any chamber $C\fl F_1$, the pair
$$([C\prec F_1],[C\prec F_2])$$
belongs to $\Phi$ if and only if the following conditions hold:\begin{enumerate}
\item[(a)] $F_2\pol F_1$, and
\item[(b)] for all $G\in\FF$ with $\codim(G)=\codim(F_1)-1$ such that $C\prec G\prec F_1$,\\ one has $G\pol F_1.$
\end{enumerate}\end{df}

We conclude by pointing out that the above definition indeed has the required features.

\begin{thm}[See Theorem 6 of \cite{salvettisette}] \label{T6}
The matching $\Phi$ is the
gradient of a combinatorial Morse function with the minimal possible number of critical cells.

\noindent Moreover, the set of $k-$dimensional critical cells is given by
$$\textup{Crit}_k(\S) = \bigg\{[C\fleq F] \left\vert\begin{array}{l} \codim(F)=k,\, F\cap V_k \neq \emptyset,\\
G\pol F \textrm{ for all }G \textrm{ with } C\fl G \precneqq F\end{array}\right.\bigg\}$$
(equivalently, $F\cap V_k$ is the maximum in polar ordering among
all facets of $C\cap V_k$).
\end{thm}


\subsection{Oriented matroids and flippings} The combinatorial data of a real arrangement of hyperplanes are customarily encoded in the corresponding {\em oriented matroid}. For the precise definition and a comprehensive introduction into the subject we refer to \cite{BLSWZ}. One of the many different ways to look at an oriented matroid is to characterize its set of {\em covectors}. Given a ground set of elements $E$, a subset of $\{-,0,+\}^{E}$ is the set of covectors of an oriented matroid if it satisfies a certain set of axioms (see \cite[Definition 3.7.5]{BLSWZ}). It is customary to partially order the set of covectors of an oriented matroid by inclusion of their support (the support of a covector $X\in\{-,0,+\}^E$ is the set of all $e\in E$ with $X(e)\neq 0$). The height of this poset (i.e., the length of every maximal chain) is the {\em rank} of the oriented matroid.

If we arbitrarily choose a positive side of every hyperplane of an arrangement $\AA$ of {\em linear} hyperplanes, we can associate to every $F\in\FF(\AA)$ the sign vector $X$ on the ground set $\AA$ with  $X(H)=+$, $-$ or $0$ if $F$ is on the positive side, on the negative side or {\em on} the hyperplane $H$. Indeed, the set of such sign vectors satisfies the axioms for the set of covectors of an oriented matroid, with the ordering of covectors naturally corresponding to the partial ordering of $\FF(\AA)$ that we defined earlier. 

However, oriented matroids are more general than linear hyperplane arrangements. To see this, recall that a $k$-pseudosphere in the $d$-sphere is the image of $S^k\subset S^d$ under a tame selfhomeomorphism of $S^d$. An arrangement of pseudospheres is a set of centrally symmetric pseudospheres arranged on the $d$-sphere in such a way that the intersection of every two pseudospheres is again a pseudosphere.  \\
The {\em topological representation theorem} (Folkman and Lawrence \cite{FoLa}, see also \cite[Theorem 5.2.1]{BLSWZ}) proves that the poset of covectors of every oriented matroid of rank $d$ can be ``represented'' by the stratification of $S^d$ induced by an arrangement of pseudospheres.
\newcommand{\TT}{\mathcal{T}}

\begin{df}[Compare Definition 7.3.4 of \cite{BLSWZ}]\label{df:flip} Let $\AA:=(S_e)_{e\in E}$ be an arrangement of pseudospheres on $S^d$. Pick a vertex $w$ of the induced stratification of $S^d$ and consider a pseudosphere $S_f$ with $w\not\in S_f$. Let $\TT_w:=\{e\in E\mid S_e\ni w\}\cup \{f\}$ and set $\UU_w:= E\setminus \TT_w$. 

We say that $w$ is {\em near} $S_f$ if all the vertices of the arrangement $\TT_w$ are inside the two regions of $\UU_w$ that contain $w$ and $-w$. 
\end{df}

Given an arrangement of pseudospheres, if a vertex $w$ is near some pseudosphere $S_f$, one can perturb locally the picture by `pushing $S_f$ across $w$' and, symmetrically, across $-w$, so to obtain another valid arrangement of pseudospheres which oriented matroid differs from the preceding only in faces inside the two regions of  $\TT_w$ that contain $w$ and $-w$. This operation was called a {\em flipping} of the oriented matroid at the vertex $w$ by Fukuda and Tamura, who first described this operation \cite{FuTa}. For a formally precise description of flippings see also \cite[p. 299 and ff.]{BLSWZ}.

Every arrangement of linear hyperplanes in $\R^d$ induces on the unit sphere $S^{d-1}$ an arrangement of spheres. An oriented matroid that can be realized in this way is called {\em realizable}. It is NP-hard to decide whether an oriented matroid is realizable \cite{JRG}. 

\begin{rem}\label{lattice} Flippings preserve the underlying matroid (i.e.,the intersection lattice of the arrangement). However, a flipping of a realizable oriented matroid need not be realizable!
\end{rem}

To be able to encode the data of an affine arrangement one uses {\em affine oriented matroids}. The idea is to add an hyperplane `at infinity' to the oriented matroid represented by the cone of the given affine arrangement (for the precise definition, see \cite[Section 4.5]{BLSWZ}). For the affine counterpart of the representation theorem we need one more definition.

\begin{df}
A $k$-pseudoflat in $\mathbb{R}^d$ is any image of $\mathbb{R}^{d-k}$ under a (tame) selfhomeomorphism of $\mathbb{R}^d$. A pseudohyperplane clearly has two well-defined {\em sides}. An arrangement of pseudohyperplanes is a set of such objects satisfying the condition that every intersection of pseudohyperplanes is again a pseudoflat. 
\end{df}

Then every affine oriented matroid is represented by an (affine) arrangement of pseudohyperplanes, and the notion of flipping is similar to the previous: the only difference is that there is no vertex ``$-w$''. 
\begin{nt}\label{nt:flip} Let $\AA$ be an affine arrangement of pseudohyperplanes, $\Ha\in \AA$, and $w$ a vertex of $\AA$ near $\Ha$. The arrangement representing the oriented matroid obtained from the previous by flipping $\Ha$ across $w$ will be denoted $\texttt{Flip}(\AA,\Ha,w)$.

\end{nt}

Consider an arrangement of affine pseudohyperplanes $\AA$ and pick a pseudohyperplane $H$ such that all points of $\AA$ are on the same side of $H$.  A {\em sweeping} (or `topological sweeping') of $H$ through $\AA$ is a sequence of flippings, one for every point of $\AA$, that fixes everything except $H$. At the end of a sweeping, the points of $\AA$ are all on the opposite side of $H$ with respect to the beginning.

It is a well-known fact that such a sweeping need not exist in general for all $\AA$ and $H$. At every step, the flip through a point $p$ of $\AA$ is performed by extending $\AA$ with a pseudohyperplane through $p$ parallel to $H$, and then perturbing the resulting arrangement around $p$ \cite[Section 7.3]{BLSWZ}. While the `perturbation' part is always feasible, the `extension' part requires careful consideration. 

The oriented matroid program $(\AA,H)$ is called {\em Euclidean} if an extension of $\AA$ by a pseudohyperplane parallel to $H$ containing $p$ exists for every point $p$ \cite[Definition 10.5.2]{BLSWZ}.  
The following characterization was first proved  in Komei Fukuda's PhD. thesis. We refer to \cite[Chapter 10]{BLSWZ} and the bibliography cited therein for a structured and complete exposition of the subject. 

\begin{tdef}[See Section 10.5, Theorem 10.5.5 of \cite{BLSWZ}]\label{teodef} Let an affine arrangement of pseudohyperplanes $\AA$ be given, and let $H\in\AA$ be such that all points of $\AA\setminus\{H\}$ are on the same side of $H$. 
Every 1-dimensional face $F$ of $\AA$ that is not contained in $H$ is supported on a pseudoline $\ell_F:=\bigcap\AA_F$, and $\ell_F$ meets $H$ in exactly one point $p$. We can then think of the 1-cell $F$ as being directed away from $p$ (along $\ell_F$). Thus, we turn the union of the $0$- and $1$- dimensional faces of $\AA$ not contained in $H$ into an oriented graph we call $G_H$. 


{\em The oriented matroid program $(\AA, H)$ is Euclidean if and only if $G_H$ is acyclic. }

\end{tdef}


\begin{crl} If an oriented matroid program $(\AA,H)$ is realizable (i.e., $\AA$ is an arrangement of hyperplanes), then $G_H$ is acyclic, and thus allows for a sweeping of $H$ through $\AA$. 
\end{crl}


\renewcommand{\thesection}{\arabic{section}}
\setcounter{section}{0}
\numberwithin{section}{part}
\numberwithin{thm}{section}
\part{Combinatorics of polar orderings}

The first step on the way to generalizing the construction of \cite{salvettisette} is to give a combinatorial (i.e., `coordinate-free') description of it. The idea is to let the hyperplane $V_{k-1}$ `sweep' across the arrangement $\AA\cap V_{k}$ instead of rotating it around $V_{k-1}$.

As explained in the introduction, we want to put the polar ordering into the broader context of the orderings that can be obtained by letting an hyperplane sweep across an affine arrangement along a sequence of flippings. By Remark \ref{lattice} we must then work with general oriented matroids, since realizability of every intermediate step is not guaranteed (and, indeed, rarely occurs). This raises the question of whether such a `sweeping' is always possible throughout the construction. We will see that indeed all occurring oriented matroid programs are Euclidean.

\newcommand{\Aa}{\widetilde{\mathcal{A}}}
\newcommand{\flip}{\texttt{Flip}}

\section{Definitions and setup}\label{sec:df}

 Let $\AA$ denote an affine real arrangement of hyperplanes in $\mathbb{R}^d$. A flag $(V_k)_{k=0,\ldots,d}$ of affine subspaces is called a {\em general flag} if every one of its subspaces is in general position with respect to $\AA$ and if, for every $k=0,\ldots d-1$, $V_k$ does not intersect any bounded chamber of the arrangement $\AA\cap V_{k+1}$. Note that this is a less restrictive hypothesis than the one required for being a {\em generic} flag in \cite{salvettisette}.

Moreover, we write 
\begin{gather*}
\AA^k:=\{H\cap V_k \mid H\in \AA \},\quad
\FF^k:=\{F\in\FF\mid F\cap V_k\neq \emptyset\}(=\FF(\AA^k)),\\
\PP^k=\{p_1,p_2,\ldots\}:=\max\FF^k,\quad \PP:=\PP^0\cup\PP^1\cup\ldots\cup\PP^d,
\end{gather*}

\noindent where of course the set $\FF^k$ is partially ordered as the face poset of the arrangement $\AA^k$.

If a total ordering $\leadsto^k$ of each $\PP^k$ is given, we define a total ordering of $\PP$ by setting, for any $p\in\PP^i$ and $q\in\PP^j$, $$p\leadsto q\Leftrightarrow \left\{\begin{array}{ll} p\leadsto^k q &\textrm{if } k=i=j \\ i<j &\textrm{if } i\neq j\end{array}\right.$$

\newcommand{\ha}{\widetilde{\HH}}

We want to let the hyperplane $V_{k-1}$ sweep across $\AA^k$. Let us introduce the necessary notation. For every $k=1,\ldots , d$, let
$$\Ha^k_0:=V_{k-1},\quad
\FF^k_0:=\FF^{k-1},\quad
\Aa^k_0:= \AA^k\cup \{\Ha^k_0\}.$$
For all $j>0$, let $p_j\in\PP^k$ be {\em near} $\Ha^k_{j-1}$ in the sense of Definition \ref{df:flip} and set
\begin{gather*}
\Aa^k_{j}:=\flip(\Aa^k_{j-1}, \Ha^k_{j-1}, p_j),\quad
\Ha^k_{j}:\quad \Aa^k_{j}\setminus \AA=\{\Ha^k_j\},\\
\HH^k_j:=(\Aa^k_j)^{\Ha^k_j},\quad \FF^k_j:=\FF(\HH^k_j),\quad
\PP^k_j:=\max \FF_j^k
\end{gather*}

\noindent where the definitions refer to  the natural inclusions $\FF^k_i\hookrightarrow \FF^k\hookrightarrow \FF$. Moreover, we will make use of the natural forgetful projection $\pi_j^k:\FF(\Aa^k_j)\rightarrow \FF^k$ (`forgetting' $\Ha^k_j$).

\begin{figure}[h]
\begin{center}
\includegraphics[scale=0.8]{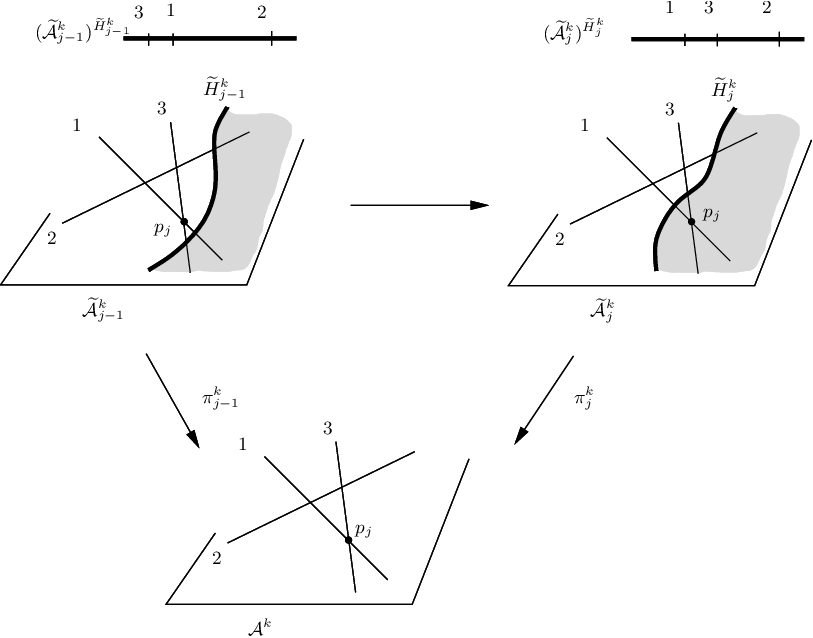}
\end{center}
\end{figure}

\begin{rem} Our construction will be inductive in the dimension. The definitions and arguments we make here about $\AA$ will be applied to every $\HH^k_j$, and so on. The involved oriented matroids can become quickly nonrealizable. Thus, it has to be stressed that our arguments hold in the generality of affine arrangements of pseudohyperplanes. The reason why we carry out this section by referring to $\AA$ as an arrangement of hyperplanes is mainly to keep the terminology lighter and help the intuition. The reader will obtain proof of the corresponding statements for pseudoarrangements by just adding throughout the next section the prefix ``pseudo'' to the  appropriate words.
\end{rem}

We have to understand how the combinatorics of the arrangement induced
on the ``moving hyperplane'' $\Ha^k_j$ changes, as $j$ becomes
bigger. By the definition of flippings, we know that
nothing changes in $\Aa^k_j$ outside
$$\YY(p_j):=(\pi^k_j)^{-1}(\FF^k_{\fleq p_j})$$\noindent - {\em a
  fortiori}, nothing changes in $\FF^k_{j-1}$ outside
$$\XX(p_j):=\FF^k_{j-1}\cap \YY(p_j).$$ 
\begin{nt} Given two faces  $F\fl G$, let us from now denote by $\op_G(F)$ the unique element of $\FF$ such that $\op_G(F)\fl G$ and the face that represents $\op_G(F)$ is on the opposite side (with respect to $F$) of every pseudohyperplane that contains $G$ but not $F$.
\end{nt}
The next Lemma states an explicit (and order-preserving) bijection between the set of `new faces'
that are cut by the moving hyperplane after the flip at $p_j$ and the following set of `old faces': $$\CC(p_j):=\{X\in \XX(p_j)\mid \op_{p_j}(X)\not\in\XX(p_j)\}.$$

\begin{lm}\label{faccedelflip} With the notations explained above, let $\Aa^k_{j-1}$ be given and let $p_j\in\PP^k$ be near $\Ha^k_{j-1}$. Then, if $<_{j-1}$ denotes the ordering of $\FF^k_{j-1}$, $\FF^k_j$ is isomorphic to the poset given on the element set $$\big(\FF^k_{j-1}\setminus \CC(p_j)\big) \cup \{(p_j,X)\mid X\in\CC(p_j)\}$$
by the order relation 
$$F\leq_{j} F^\ast :\Leftrightarrow\left\{\begin{array}{ll}
F,F^\ast\in\FF^k_{j-1}\setminus\CC(p_j) & \textrm{and } F\leq_{j-1} F^\ast,\\
F=(p_j,X), F^\ast=(p_j,X^\ast)& \textrm{and }X\leq_{j-1} X^\ast,\\
F=(p_j,X), F^\ast\in\FF^k_{j-1}\setminus\CC(p_j)& \textrm{and }\op_{p_j}(X)\leq_{j-1} F^\ast,
\end{array}\right.$$
the isomorphism being given by the correspondence $(p_j,X)\mapsto \op_{p_j}(X)$, and the identical mapping elsewhere.
\end{lm}

\begin{pf} Compare  \cite[Corollary 7.3.6]{BLSWZ}.
\qed\end{pf}

Note that the faces represented by $(p_j,X)$ for $X\in\CC(p_j)$ are exactly the faces $F$ whose minimal $k$-face is $p_j$.


\begin{crl}\label{daspostare} If $p_i,p_{i+1}\in \PP^k$ are both near $\Ha_{i-1}^k$, then the structure of $\Aa_{i+1}^k$ does not depend on the order in which the two flippings are carried out.  

In particular, any $q\in \PP^k$ near $\Ha_{i-1}^k$ and different from $p_i$ is also near $\Ha_i^k$.
\end{crl}
\begin{pf} The fact that both are near $\Ha_{i-1}^k$ implies in particular $\CC(p_i)\cap \CC(p_j)=\emptyset$, and thus the modifications do not influence each other.
\qed\end{pf}

\begin{nt}

Every $\HH^k_j$ contains an isomorphic copy of $\FF^{k-1}_0\simeq \FF^{k-2}$
because $\FF(\HH^k_0)=\FF^{k-1}$. We may then add to $\HH^k_j$ a pseudohyperplane $\lkj_0$ that intersect exactly the faces of $\FF^{k-2}$ (`a copy of $\FF(\HH^{k-1}_0)$') and consider consecutive flippings $\lkj_i$ of it along the elements of $\PP^k_j$. 
\end{nt}

\begin{rem}\label{OMP_expl} It is not difficult to see that $\La^{k,j}_0$ indeed {\em can} be swept through $\HH_j^k$. First of all, the oriented matroid program defined by $\HH_0^{k}$
and $\La^{k,0}_0$ is euclidean because the oriented matroid associated to $\HH^k_0$
is realizable (this arrangement is obtained by intersecting $V_{k-1}$ with $\AA$).  To conclude that $\lkj_0$ can be swept through $\HH_j^k$ for $j>0$ it is enough to see that, for every $j\geq 0$, euclideanness of the program associated with $\HH_j^k$ and $\lkj_0$ implies euclideanness of the program associated with $\HH_{j+1}^k$ and $\La^{k,j+1}_0$.

This last fact is readily checked by considering in both cases the orientation of the graph associated to the programs. By Lemma \ref{faccedelflip} we know how $\HH_j^k$ changes to $\HH_{j+1}^k$ after the flip through $p_j$, and since $\La^{k,j}_0=\La^{k,j+1}_0$, the orientation of the edges agrees everywhere except in $\CC(p_j)$. Now by inspecion of the possible situations one concludes that the existence of a directed cycle in the graph associated to $\HH_{j+1}^k,\,\La^{k,j+1}_0$, implies the existence of a directed cycle in the graph associated to $\HH_{j}^k,\, \La_0^{k,j}$. Then, by \ref{teodef} we are done.
 \end{rem}

\section{Special orderings}\label{special}

\begin{df} Given an essential affine real (pseudo)arrangement $\AA$ and a general
  position (pseudo)hyperplane $\Ha_0$, a total ordering $p_1,p_2,\ldots$ of
  the points of $\AA$ is a {\em special ordering} if there is a
  sequence of arrangements of pseudohyperplanes $\Aa_0,\Aa_1,\ldots$
  such that
  $\Aa_0=\AA\cup\{\Ha_0\}$, and for all $j>0$, $\Aa_j$ is obtained
  from $\Aa_{j-1}$ by flipping $\Ha_j$ across $p_j$. 
\end{df}

We collect some fact for later reference.

\begin{rem}\label{duecond} It is clear that every $\Ha^k_j$ is in general position with respect to $\AA$, because $\Ha^k_0$ was chosen so. Therefore, any two $p
,q$ that are near some $\Ha^k_j$ satisfy $\CC(p)\cap \CC(q)=\emptyset$ (just by definition of `near', see \cite{BLSWZ}).  This means amongst other that every element of $\FF_{\preceq p}\cap \FF_{\preceq q}$ is already in $\HH^k_{j}$, thus either is in $V_{k-1}$ or in some `earlier' $\CC(z)$, for $z\leadsto^k p_j\leadsto^k p,q$.
 \end{rem}

\begin{lm}\label{tripletta}
Let a special ordering $\leadsto$ of the points of an affine arrangement $\AA$ with respect to a generic hyperplane $\Ha_0$ be given. Choose two consecutive points $p\leadsto q$ and let $\leadsto^\ast$ be the total ordering of obtained from $\leadsto$ by reversing the order of $p$ and $q$. Then, the following are equivalent:\begin{enumerate}
\item[(1)] $\leadsto^\ast$ is a special ordering with respect to $\Ha_0$,
\item[(2)] In the induced flipping sequence just before the flipping through $p$, both $p$ and $q$ are near the moving pseudohyperplane. 
\item[(3)] For all $F\in\FF_{\fleq p}\cap\FF_{\fleq q}$, the minimum vertex of $F$ comes before $p$ and $q$ in $\leadsto$. 
\end{enumerate} 
\end{lm}
\begin{pf} (1)$\Leftrightarrow$(2) is clear, and (2)$\Leftrightarrow$(3) follows from Remark \ref{duecond} above.
\qed\end{pf}

Let us return to the setup of Section \ref{sec:df} and fix $k\in\{1,\ldots,d\}$ for this section. We want to understand whether (and how) it is possible to deduce a valid special ordering of the elements of $\PP^k_{j}$ from  a special ordering of the elements of $\PP^k_{j-1}$.

\begin{df}\label{indotto}
Let a total ordering $\leadsto^k_{j-1}$ of $\PP^k_{j-1}$ be given. For
every line $\ell$ of $\HH^k_{j-1}$ that contains some element of
$\XX(p_j)\cap\PP^k_{j-1}$ let $y^+(\ell)$, $y^-(\ell)$ denote the
points of $\HH^k_{j-1}$ where $\ell$ intersects the (topological)
boundary of $\XX(p_j)$, ordered so that $y^+(\ell)\leadsto^k_{j-1}
y^-(\ell)$.\\
Moreover, call $\overline{y}$ the maximum with respect to
$\leadsto^k_{j-1}$ of all $y^+(\ell)$ (for varying $\ell$).

Then define a total ordering of $\PP^k_{j}$ by setting, for every
$z_1,z_2\in\PP^{k}_{j}$:
$$z_1\leadsto^k_j z_2 \Leftrightarrow \left\{\begin{array}{ll}
z_1,z_2\in\PP^k_j\cap\PP^k_{j-1} & \textrm{and }z_1\leadsto^k_{j-1}z_2\\
z_1\not\in\PP^k_{j-1},z_2\in\PP^k_{j-1}& \textrm{and } \overline{y}
\leadsto^k_{j-1} z_2\\
z_i=(p_j,x_i)\textrm{ for }i=1,2&\textrm{and }x_2^\ast\leadsto^{k-1}x_1^\ast,
\end{array}\right.
$$
where $x_i^\ast$ denotes the unique element of $\PP^{k-1}$ with the same support
as $x_i$.
\end{df}

Our goal will be to prove the following statement.

\begin{thm}\label{crescej} For every $k\geq 0$ and every $j>0$, if $\leadsto^k_{j-1}$ is a special ordering, so is $\leadsto^k_{j}$ too. 
\end{thm}

\begin{nt}To investigate the situation, we will focus on
  $\XX(p_j)\subset\HH^k_{j-1}$.  Let us write $x_1,\ldots, x_s$ for the points of this complex. Also, let $\ell_1,\ldots,\ell_l$ be the (pseudo)lines of $\HH^k_j$ that contain some $x_i$ and write $y_1,y_2,\ldots$ for the intersection points of the $\ell$'s with the hyperplanes bounding $\XX(p_j)$.\end{nt}

\begin{figure}[h]
\begin{center}
\includegraphics[scale=0.7]{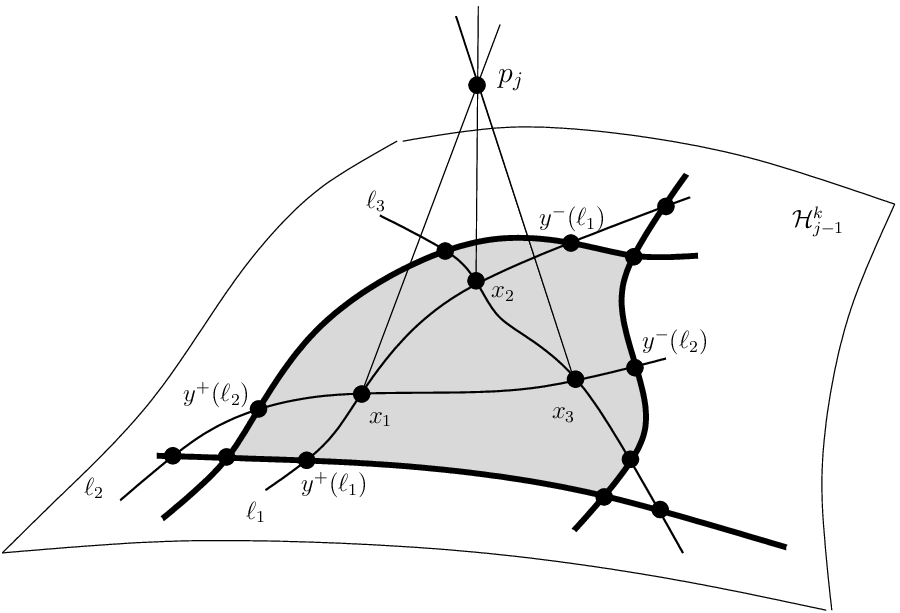}
\caption{An illustration of our setup. The shaded region is $\XX(p_j)$, and the subcomplex $\CC(p_j)$ is spanned by $x_1,x_2,x_3$.}\end{center}
\end{figure}

\begin{rem}\label{inriga} It is useful to consider
  the { lines} passing through a point $q\in\PP^k$. For instance,
  one can see that if two points $p,q\in\PP^k$ lie on a common line
  $\ell$ of $\AA^k$ so that $p$ is nearer than $q$ to $\ell\cap V_{k-1}$, then there is no sequence of flippings of $\Ha^k_0$ in which $q$ comes before $p$.
\end{rem}

\begin{lm}\label{isolati} Let a special ordering of $\PP^k_{j-1}$ be
  given. Let $\XX(p_j)=\{x_1,\ldots, x_s\}$ be numbered so that  $V_{k-1}\cap \vert
  x_r\vert\leadsto^{k-1} V_{k-1}\cap \vert x_t\vert$ if and only if $r<t$ (remember that $\vert x\vert$ denotes the support of $x$). 
Moreover, let $p_1,p_2,\ldots$ denote the elements of $\PP^k_{j-1}\setminus\{x_1,\ldots, x_s\}$ ordered according to $\leadsto^{k}_{j-1}$ and let $m$ be so that $p_m=\overline{y}$. 
Then the following is a special ordering of $\PP^k_{j-1}$:
$$p_1,p_2,\ldots, \overline{y},x_1,x_2,\ldots,x_s,p_{m+1},p_{m+2},\ldots.$$

\end{lm}

\begin{pf} The proof is subdivided in three steps.

\begin{claim}\label{banale} Every $y_i$ is contained in exactly one of the lines $\ell_1,\ldots, \ell_l$. Moreover,  for all $1\leq i<j\leq l$, there is $r$, $1\leq r\leq s$, such that $x_r=\ell_i\cap\ell_j$.
\end{claim}
\begin{pf1}
Note that $\ell_i\cap\ell_j\neq \emptyset$ because both lines are flats of the central arrangement $\AA_{p_j}$, and these intersections are points of the arrangement $\Ha^k_{j-1}\cup\AA_{p_j}$. Now both claims follow because the subcomplex $\XX(p_j)$ contains, by definition of flipping, every point of the arrangement given by $\Ha^k_{j-1} \cup \AA_{p_j}$ (see Definition \ref{df:flip} and ff.).
\end{pf1}

Now recall that,  in {\em any} special ordering of $\PP^k_{j-1}$, the $0$-dimensional faces on every $\ell_i$ must be ordered `along $\ell_i$'. Thus, on every line $\ell_i$ the segment contained in $\XX(p_j)$ is bounded by two points, say $y^+(\ell_i) \leadsto^k_{j-1} y^-(\ell_i)$.

 \begin{claim}\label{claim2} Consider a special ordering of $\PP^k_{j-1}$. Then the ordering remains special after the following modifications:
\begin{enumerate}
\item[(1)] Switching $y^+(\ell)$ and $x$ whenever $x$ comes right before $y^+(\ell)$.
\item[(2)] Switching $y^-(\ell)$ and $x$ whenever $x$ comes right after $y^-(\ell)$.
\item[(3)] Switching $x$ and any $z\not\in\XX(q)$ whenever $x$ and $z$ are consecutive.
\end{enumerate}\end{claim}

\begin{pf2}
In case (1) note that Claim \ref{banale} ensures that $\CC(y^+(\ell))$
lies fully outside $\XX(p_j)$ and so it is disjoint from any
$\CC(x)$. Now let $x$ be, say, the $r$-th element of $\PP^k_{j-1}$. Since $x$ comes
right before $y^+(\ell)$ we must have that $y^+(\ell)$ is already near
$\La^{k,j-1}_{r-1}$: indeed, in that case $x$ cannot be contained in
$\ell$ and by definition also not in the boundary hyperplane that
intersects $\ell$ in $y^+(\ell)$. Since the only change in passing
from $\La^{k,j-1}_{r-1}$ to $\La^{k,j-1}_{r}$ happens at faces which
supports contain $x$, we have $\YY(y^+(\ell))\cap \La^{k,j-1}_{r-1}=\YY(y^+(\ell))\cap\La^{k,j-1}_{r}$. By Corollary \ref{daspostare} we are done.

The case (2) is handled similarly, by reversing the order of the flippings, and case (3) is clear. \qed\end{pf2}

At this point we know that the ordering
$$p_1,p_2,\ldots,p_m,[\cdots ],p_{m+1},p_{m+2},\ldots,$$
where the square brackets contain the $x_i$'s, is indeed a special ordering of $\PP^k_{j-1}$. We have to prove that we can indeed arrange the elements in the square bracket as required. 

First, if $x_1$ is not near  $\La^{k,j-1}_m$, then there is a line 
$\ell\ni x_1$ and some other $x_i$ that lies on $\ell$ between $x_1$ and 
$\ell\cap\La^{k,j-1}_m$. In particular, $x_i$ lies between $x_1$ and 
$\ell\cap\La^{k,j-1}_0=\ell\cap\FF^{k-1}_0=\ell\cap V_{k-2}$. 
The points $x_1,\ldots,x_s$ are given by the intersection of the 
pseudohyperplane $\HH^k_{j-1}$ with lines $g_1,\ldots,g_s$ of $\AA^k$, 
and $\ell$ is the intersection of $\HH^k_{j-1}$ with the plane $E$ generated by
 $g_1$ and $g_i$. For all $r$ let $x_r^\ast:=g_r\cap V_{k-1}$. 
Since $g_1\cap g_i=p_j$, that lies outside the segments 
$\overline{x_1 x_1^\ast}$ and $\overline{x_i x_i^\ast}$, we get that in 
$V_{k-1}$ the point $x_i^\ast$ lies on the line $\ell^\ast:=E\cap V_{k-1}$ 
between $x_1^\ast$ and $\ell^\ast\cap \Ha^{k-1}_0= \ell^\ast \cap V_{k-2}$. 
With Remark \ref{inriga}, and by the way the numbering of the $x_r$ was chosen,
we reach a contradiction. We may now repete the argument with $x_2$, and all the following points until we reach $x_s$, concluding the proof.
\end{pf}


We are now ready to prove the main result of this section.

\begin{pft} We can assume that
  $\leadsto^k_{j-1}$ is modified so to agree with the statement of
  Lemma \ref{isolati}. Let $U^{k,j}_m:=\bigcup_{i\leq m}
  \La^{k,j}_i$ (meaning the set of all faces that are contained in some $\La^{k,j}_i$). Since the orderings $\leadsto^k_{j-1}$ and
  $\leadsto^{k}_{j}$ now agree up to $p_m=\overline{y}$ and clearly
  $U^{k,j}_m= U^{k,j-1}_m$ by Lemma \ref{faccedelflip}, we are left
  with proving that it is possible to perform the flippings of the
  $x_i$ just after $\overline{y}$, and in the reverse order as the
  corresponding flippings are performed in $\Ha^k_{j-1}$.

To this end, let us consider $\La^{k,j}_m$, i.e., the moving pseudohyperplane `just after' the flipping through $p_m=\overline{y}$. Recall that $\La^{k,j}_m\simeq \La^{k,j-1}_m$, and in particular we can compare the points $z_1,\ldots, z_l$ where the lines containing some $x_i$ intersect the pseudohyperplane corresponding to $\La^{k,j}_m$. Let $F_1,\ldots,F_l$ be the faces such that $z_i=F_i\cap\La^{k,j-1}_m$. Then we see that the `same' points $z_i$  are given by $(p_j,F_i)\cap\La^{k,j}_m$. So by the correspondence established in Lemma \ref{faccedelflip} we have that a point $(p_j,x)$  is near $\La^{k,j}_m$ if and only if $x$ is near (but ``on the backside'' of) $\La^{k,j}_{m+s}$. This shows that $(p_j,x_s)$ is near $\La^{k,j}_m$. After performing this flipping we may repeat the argument to conclude that $(p_j,x_{s-l})$ is near $\La^{k,j}_{m+l}$ for every $l\leq s$, and the claim of the Theorem follows. \end{pft}

\section{Combinatorial polar orderings}\label{combinatorial}

After having looked inside each $V_k$, let us study the structure that arises by considering all strata. 


\begin{df}[Compare Theorem 5. of \cite{salvettisette}]\label{cinque} Given 
total orderings $\leadsto^k$ of each $\PP^k$, we define a total ordering 
$\pol$ of $\FF$. All faces of codimension $d$ are elements of $\PP^d$ and are ordered accordingly. Assuming the ordering is defined for all faces of codimension  $k+1$ and bigger, then given two $k$-codimensional faces $F$ and $G$ we have:
\begin{enumerate}
\item[(1)] if $F,G\in\PP^k$, $F\pol G$ if $F\leadsto G$;
\item[(2)] if $F\in\PP^k$ and $G \not\in\PP^k$, then $F\pol G$; 
\item[(3)] if $F,G\not\in\PP^k$, let $F'$, (resp. $G'$) be 
the $k+1$-codimensional facet in the boundary of $F$ (resp. $G$), which is 
minimum with respect to $\pol$. Then: 
\begin{enumerate}
\item[(3.1)] if $F'\pol G'$, then $F\pol G$.
\item[(3.2)] if $F'=G'$, then $F\pol G$ if and only if 
$F_0\leadsto G_0$, where $F_0$ and $G_0$ are the unique elements of 
$\PP^{k}$ that have the same linear span as $F$, respectively $G$. 
\end{enumerate}
\item[(4)] If $F\in\PP^k$, then F is lower than any $k+1$-codimensional facet 
\item[(5)] If $F\not\in\PP^k$, then $F$ is bigger than its minimal boundary
$F'$ and lower than any $(k+1)$-codimensional facet which is bigger 
than $F'$.
\end{enumerate}
\end{df}

Thus, if the orderings on the $\PP_k$s are given by lexicografic order on the 
polar coordinates, we reproduce the polar order of \cite{salvettisette}.\\

\begin{df} Let an affine real arrangement $\AA$ be given. A {\em combinatorial polar ordering} of $\FF(\AA)$ is any total ordering $\pol$ induced via Definition \ref{cinque} by the choice of a general flag $(V_k)_{k=0\ldots d}$ and of special orderings $\leadsto^k$ of the points of $V_k$ with respect to $V_{k-1}$, for every $k=1,\ldots, d$. 
\end{df}

Let us next give an alternative characterization of the combinatorial polar orderings that will turn out to be useful later on.

\begin{df}
Given $F\in\FF$, define the signature of $F$ as
$\sigma(F)=(k_F,j_F,m_F)$, where
\begin{gather*} 
k_F:=\min\{k\mid V_k\cap F\neq \emptyset \}\\
j_F:=\min\{j\mid F\in\FF(\HH^{k_F}_j)\}\\
m_F:=\min\{m\mid F\in \FF(\La^{k_F,j_F}_m)\},
\end{gather*} 
\noindent where we agree to put $j_F=0$ when $k_F=0$ and $m_F=0$ if $k_F\leq 1$ because in those cases the above definition is void.
\end{df}

\begin{lm}\label{trilex}
Let special orderings $\leadsto^k$ be given for every $k$, and let $\pol$ be the total ordering of $\FF$ induced by them. For $F_1,F_2\in\FF$, if $\sigma(F_1)<\sigma(F_2)$ in the lexicographic order, then $F_1\pol F_2$.
\end{lm}
\begin{pf} If $k_{F_1}<k_{F_2}$, then by Definition \ref{cinque}.(4) $F_1\pol F_2$. 

Suppose now $k_{F_1}=k_{F_2}$ but $j_{F_1}< j_{F_2}$. If $F_1,F_2\in\PP^k$, then we are already done by Definition \ref{cinque}.(1). Else, the condition means that the minimal codimensional-$k+1$ face of $F_1$ comes before the minimal codimensional-$(k+1)$ face of $F_2$, and by Remark \ref{inriga} we are done.

The same line of reasoning applies to show that $k_{F_1}=k_{F_2}$, $j_{F_1}=j_{F_2}$ and $m_{F_1}<m_{F_2}$ implies $F_1\pol F_2$. \qed\end{pf}

\begin{rem}\label{lex} It is now easy to see that one could go on and
  define for every face $F$ a vector $$(\sigma_1(F),\ldots,
  \sigma_{k_F}(F))$$ with $\sigma_1(F):=j_F$ and
  $\sigma_i(F):=\min\{m\mid F\in
  \La_m^{k_F,\sigma_1(F),\ldots,\sigma_{i-1}(F)}\}$ (where
   $\La_m^{k_F,a_1,a_2,\ldots,a_j}$ is defined for $j>1$ as the moving hyperplane of
  $\HH^{k_F,a_1,\ldots,a_{j-1}}_{a_{j}}$ after the $m$-th flipping). From this, a signature 
$$\sigma(F):=(\underbrace{0,0,\ldots,0}_{d-k_F \textrm{ times}},\sigma_1(F),\ldots,\sigma_{k_F}(F))$$
can be defined, so that for all $F_1,F_2\in\FF$, $F_1\pol F_2$ if and
only if $\sigma(F_1)<\sigma(F_2)$ lexicographically. This yields an
alternative equivalent formulation of the ordering defined in \ref{cinque}.
\end{rem}

\begin{rem} From the point of view of the computational complexity,
  the translation of Remark \ref{lex} shows that the whole work
  amounts indeed to determine special orderings of the
  $V_k$'s. Effective algorithms for this kind of tasks were developed
  in the last few years by Edelsbrunner et al. \cite{Edels}.
\end{rem}

\section{``Polar'' vector fields and switches}\label{polar}

Recall that for $F\in\FF$ we denote by $F'$ the smallest facet of $F$ with respect to the given ordering $\pol$. We rephrase Definition \ref{T6} in our broader context.

\begin{df}\label{Fi} Let an affine real arrangement $\AA$ and a general flag
  $\{V_k\}_{k=0,\ldots,d}$ be given. For every total ordering $\pol$
  of $\FF$ we define
$$\Phi(\pol):=\left\{\begin{array}{crl} & (i)&\!\!\!\! F\not \in \PP,\\
 {[C\fleq F]< [C \fleq F']}\in\S :& (ii)& \!\!\!\! G'\neq F
\textrm{ for all } G
\textrm{ with }\\&& \!\!\!\! C\fl G\fl F.\end{array}\right\}.$$
\end{df}

\begin{rem} If $\pol$ is the polar ordering defined in \cite{salvettisette}, then
  by Theorem \ref{T6} we know that $\Phi(\pol)$ is a
  maximum acyclic matching on the poset of cells of the Salvetti
  complex, i.e., it defines a discrete Morse function on $\S$ with
  the minimum possible number of critical cells. 
\end{rem}

Our aim is to show that the total ordering can be slightly modified
without affecting the resulting acyclic matching.

\begin{df}[Switch]\label{switch}
Let special orderings $\leadsto^k$ of the $\PP^k$'s with respect to
$V_{k-1}$ be given and let $\pol$ denote the induced total ordering of $\FF$. \\
Two faces $F_1 F_2\in\PP^k$ are called {\em c-independent} if

(1) they are consecutive with respect to $\leadsto^k$, and

(2) $G\pol F_1,F_2$ for every $G\in\FF_{\fleq F_1}\cap\FF_{\fleq F_2}$.

\noindent The ordering $\leadsto^*$ is obtained from $\leadsto$ by a {\em switch} if there are two c-independent faces $F_1\leadsto F_2$ so that $F_2\leadsto^* F_1$, while $F\leadsto G$ implies $F\leadsto^* G$ for every other $F,G$. We will write $\pol^\ast$ for the corresponding combinatorial polar ordering.
\end{df}

The following fact is an easy consequence of Corollary \ref{daspostare}.

\begin{thm} If an ordering $\leadsto$ of the points of an affine
  arrangement is special with respect to a general position hyperplane
  $\Ha$, then so is $\leadsto^*$.
\end{thm}

Now we need to study how the induced total orderings $\pol$ of $\FF$
vary by switching two c-independent faces.

\begin{lm} Let a special ordering $\leadsto$ of $\PP$ be given, and $\pol$ be the associated total ordering of $\FF$. Moreover, let $\leadsto^*$ be obtained from $\leadsto$ by a switch and let $\pol^*$ be defined accordingly. Then the minimum facet $F'$ of any $F\in\FF$ with respect to $\pol$ is also the minimum facet with respect to $\pol^*$.
\end{lm}
\begin{pf} Let $F_1,F_2$ denote the two faces involved in the switch, and write $k_0:=k_{F_1}=k_{F_2}$. The claim is easily seen to be true if $k_F<k_0$ or if $k_F>k_0+1$.

Consider the case where $k_F=k_0$. Since the ordering $\leadsto^{k_0-1}$
does not change, if
\begin{equation}\label{minv}\min_{\leadsto}\{p\in\PP^{k_0}\mid p\fgeq
  F\} = \min_{\leadsto^*}\{p\in\PP^{k_0}\mid p\fgeq F\}\end{equation}
then the claim is clearly true by Lemma \ref{trilex}.

Because $F_1,F_2$ are consecutive, condition (\ref{minv}) fails only if both $F_1,F_2 \fg F$. But then by Definition \ref{switch}.(2) $F\pol F_1,F_2$, implying that the minimum facet of $F$ comes before $F_1$ and $F_2$, and thus remains unchanged by passing from $\pol$ to $\pol^\ast$.

Now let $k_F=k_0+1$. If $\codim(F)=k_0$, then $F'$ (i.e., the minimal facet of $F$) is an element of $\PP^{k_0+1}$,
where the order remains unchanged; in any other case,
$j_{F'}=j_{F}$. So after Lemma \ref{trilex} we must prove that the
claim holds for $F\in\op_{p_j}\CC(p_j)$, for any
$p_j\in\PP^{k_0+1}$. Because the $F_i$ are consecutive, the ordering
on the set $\PP^{k_0+1}_{j-1}\cap\XX(p_j)$ does not change in passing
from $\leadsto$ to $\leadsto^*$, unless $p_j$ is the intersection
point of the two lines of $\AA^{k_0+1}$ that contain $F_1$ and
$F_2$. But even in this last case, the corresponding points $G_1,G_2$
of $\HH^k_j$ are again consecutive. Moreover, they are not joined by
an edge in $\HH^k_j$because $F_1$ and $F_2$ are not.  By the
construction of Lemma \ref{isolati}, all this implies that they are
both near the moving pseudohyperplane $\La^{k_F,j}$ `just before
flipping across the first of them'. In turn, this means (by Remark
\ref{duecond}) that the elements of $\FF_{\fleq G_1}\cap\FF_{\fleq
  G_2}$, and in particular $F$ and $F'$, come before $G_1$ and $G_2$ - i.e., the only elements of $\PP^{k_F}_j$ that are switched. We can then apply the same reasoning as the case $k_0=k_F$ to conclude the proof. \qed\end{pf}

In particular, just by looking at the definition of the matchings we obtain the following result.

\begin{thm} Let a special ordering $\leadsto$ of $\PP$ be given, and $\pol$ be the associated total ordering of $\FF$. Moreover, let $\leadsto^*$ be obtained from $\leadsto$ by a switch and let $\pol^*$ be defined accordingly. Then
$$\Phi(\pol)=\Phi(\pol^*).$$
\end{thm}

The next step is to see that actually switches are rather powerful
tools for transforming special orderings.

\begin{thm} Let $\leadsto_1, \leadsto_2$ be any two special orderings of the point of an arrangement $\AA$ with respect to a generic hyperplane $\Ha$. Then $\leadsto_2$ can be obtained from $\leadsto_1$ by a sequence of switches.
\end{thm}

\begin{pf} Let $\PP$ denote the set of points of $\AA$. Write $\PP=\{p_1,p_2,\ldots,p_m\}$
where $i<j$ if $p_i\leadsto_1 p_j$. 
Let $\sigma$ be the permutation of $[m]$ so that $p_i\leadsto_2 p_j$ if $\sigma(i)<\sigma(j)$. We proceed by induction in the number $u(\sigma)$ of inversions in $\sigma$, the case $u(\sigma)=0$ being trivial.

So suppose $u(\sigma)>0$. Then there are numbers $i_1<i_2$ such that $\sigma(i_1)=\sigma(i_2)+1$.  If $\tau$ is the transposition $(\sigma(i_2),\sigma(i_1))$, then the number of inversions of the permutation $\tau\sigma$ is strictly smaller than $u(\sigma)$. 

Clearly the ordering of $\PP$ associated to $\tau\sigma$ is obtained by changing the position of $v_1:=p'_{\sigma(i_1)}$ and $v_2:=p'_{\sigma(i_2)}$. Thus we will be done by showing that this is a valid `switch' in $\leadsto_2$ according to Definition \ref{switch}. 

To this end, first remark that the elements are clearly consecutive in $\leadsto_2$. Next consider the fact that $v_2\leadsto_1 v_1$ and $v_1\leadsto_2 v_2$, where both $\leadsto_1$ and $\leadsto_2$ are valid special orderings.  By Remark \ref{inriga} there is no line containing both $v_1$ and $v_2$. Thus, in the sequence of flippings associated to $\leadsto_2$, just before flipping across $v_1$ the moving hyperplane is actually also near $v_2$. By Lemma \ref{tripletta} this ensures condition (2) of the definition of independence, and concludes the proof. \qed\end{pf}

If $\pol$ is the polar ordering defined in \cite{salvettisette}, 
then by Theorem~\ref{T6}. we know that $\Phi(\pol)$ is a
  maximum acyclic matching on the poset of cells of the Salvetti
  complex, i.e., it defines a discrete Morse function on $\S$ with
  the minimum possible number of critical cells. Moreover, the
  critical cells are given in terms of $\pol$ by Theorem~\ref{T6}. 

At this point, the main result of this section is evident.

\begin{prop}\label{aaa} Let a combinatorial polar ordering of the faces of an affine real arrangement $\AA$ be given. Then the induced matching $\Phi(\pol)$ is a discrete Morse vector field with the minimum possible number of critical cells. 
\end{prop}

\begin{rem} We already saw that the approach via flippings makes it unnecessary to request the stronger form of `generality' for the flag $(V_k)_k$ that is needed in \cite{salvettisette}. However, if this condition {\em is} satisfied, then the matching is the polar gradient of \cite{salvettisette}.
\end{rem}

\part{Recursively orderable arrangements}

Having established that every special ordering of an arrangement with respect to a general flag gives rise to a combinatorial polar ordering - and thus to a minimal model for the complement of the arrangement's complexification, the problem of actually finding such an ordering remains.

However, some arrangements admit some particularly handy special orderings, that give rise to combinatorial polar ordering that appear particularly well-suited for explicit computations. The motivating example here is the braid arrangement,  studied in \cite{salvettisette}. In the following we state this nice property  and look for other examples of arrangements that enjoy it.

\section{The definition}\label{the}

\newcommand{\totl}{\sqsubset}

\begin{df}[Recursive Ordering]\label{df:maindef} Let $\AA$ be a real arrangement and $(V_k)_{k=0,\ldots, d}$ a general flag. The corresponding {\em recursive ordering} is the total ordering $\totl$ of $\PP$ given by setting $F\totl G$ if one of the following occurs:\begin{enumerate}
\item[(i)] $F\in\PP^h$, $G\in\PP^k$ for $h<k$.
\item[(ii)] there is $k$ so that $F,G\in\PP^k$ and, writing $F_0:=\min \{J\in\PP^{k-1}\mid F\subset \vert J\vert\}$,  $G_0:=\min \{J\in\PP^{k-1}\mid G\subset \vert J\vert\}$, \begin{enumerate}
\item[(a)] either $F_0 \totl G_0$,
\item[(b)] or $F_0=G_0$ and there exists a sequence of faces
$$F_0 \fl F_1 \fg J_1 \fl F_2 \fg J_2 \cdots \fl F$$

 such that $\codim(F_i)=\codim(J_i)+1=\codim(F)$, and every $J_i$, $F_i$ intersect $|F_0| \cap V_k$, and $F_i\neq G$ for  all $i$.
\end{enumerate}\end{enumerate}
\end{df}

\begin{df} An arrangement $\mathcal A$ in $\R^n$ is said to be 
{\em recursively orderable} if there is a general flag $(V_k)_{k=0,\ldots,d}$ so that the corresponding recursive ordering is special.
\end{df}


\begin{ex} The braid arrangement on $n$ strands is recursively orderable for every $n$, as was shown (and exploited) in \cite{salvettisette}.
\end{ex}

\begin{rem}\label{trovarli} With the work done so far, we see that proving that an arrangement $\AA$ is recursively orderable amounts essentially to finding a special ordering of $\PP(\AA)$ such that in every $V_k$ condition (ii).(a) of the above Definition \ref{df:maindef} holds, since Conditions (i) and (ii).(b) are ``standard features'' in every special ordering. 
\end{rem}

\section{Recursively orderable arrangements of lines}\label{follow}

In this section $\AA$ will be an affine arrangement of lines in $\mathbb{R}^2$. And we will suppose it to be {\em actually} affine, i.e. $\PP^2$ consists of more than one element (otherwise the arrangement is central, and every central $2$-arrangement is trivially recursively orderable). Here we do not need the detailed notation of the general case, so we will write $P:=\PP^2$ and abuse notation by writing $\AA:=\PP^1$.

The generic flag here is a pair $(b,\ell)$, where $b$ is a point in an unbounded chamber and $\ell\ni b$ is a line in general position with respect to $\AA$ where all the points of $\AA$ lie on the same side of $\ell$, and the points $\AA\cap\ell$ lie on the same halfline with respect to $b$. We shall sometimes confuse $b$ with the chamber $B$ it is contained in. In particular, we see that $B$ cannot have two parallel walls.

\begin{nt} Let an affine arrangement of lines $\AA$ be given together
  with a general flag $(b,\ell)$. 
The line
  $\ell$ intersects a facet of $B$: let $h_0$ denote the element
  of $\AA$ supporting it. Let $a_1,a_2,\ldots$ denote the
  points on $h_0$, numbered by increasing distance from
  $b$. Moreover, write $M_j:=\{h^j_1,h^j_2,\ldots,h^j_{\max}\}$ for the set of 
  all lines different from $h_0$ that contain $a_j$, ordered according to
  the sequence of points they generate on $\ell$. For every $h\in\AA$
  let $h^+$ denote the (open) halfplane bounded by $h$ and containing
  $b$, and set $h^-:=\mathbb{R}^2\setminus h^+$. Then we define, for
  every $j=1,\ldots r$,  

\begin{gather*}\Lambda_1:= \overline h_0^+\cap (h^1_{\max})^-,\\
  \Lambda_j:=(h^{j-1}_{\max})^+\cap (h^j_{\max})^- \textrm{ for }j>1,\end{gather*}
\noindent where overline denotes topological closure.\end{nt}
\begin{figure}[h]
\begin{center}
\includegraphics[scale=0.6]{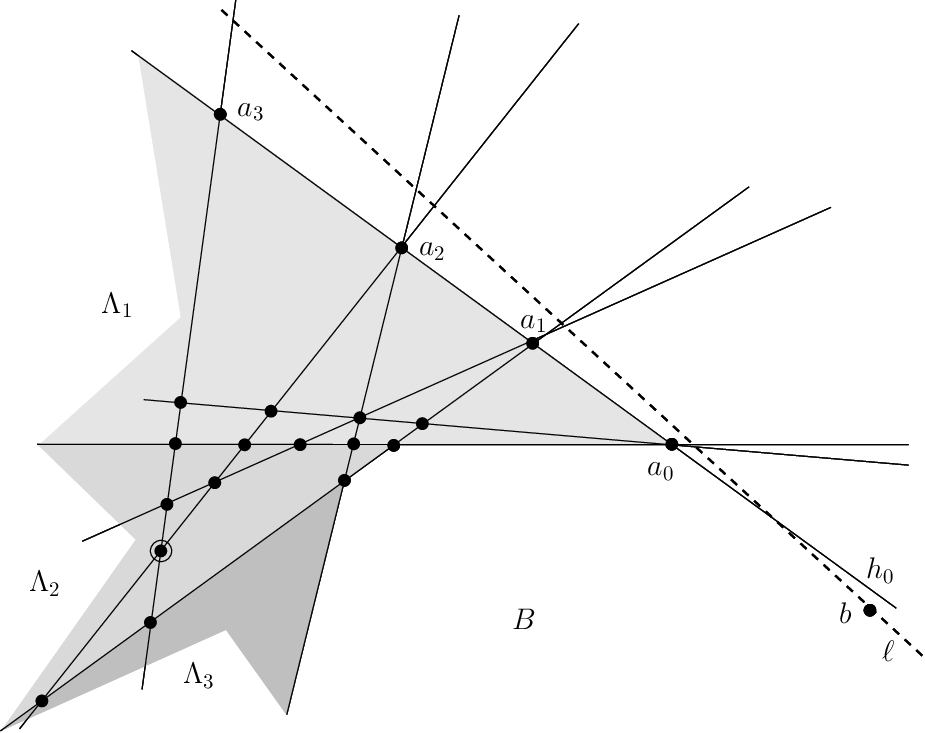}
\caption{An affine line arrangement where $\Lambda_1$ is complete with respect to $(b,\ell)$ but $\Lambda_2$ is not. Thus, it is not recursively orderable}
\end{center}
\end{figure}

\begin{df} If for every $p\in P\cap \Lambda_j$ there is $h\in M_j$ with $a_j,p\in H$, then we will say that $\Lambda_j$ is {\em complete} (with respect to $(b,\ell)$). The arrangement $\AA$ is {\em complete with respect to $(b,\ell)$} if every $\Lambda_j$ is complete and $P\subset\bigcup_{j=1,\ldots,r} \Lambda_j$
\end{df}

\begin{lm} An affine line arrangement $\AA$
  is recursively orderable with respect to a general flag $(b,\ell)$ if and
  only if $\AA$ is complete with respect to $(b,\ell)$.
\end{lm}

\begin{pfsk}
Fix an $\ell$. If $\AA$ is not complete at some $j$, then there is a point $x\in P$
 so that $x\in\Lambda_j$ but there is no line containing $a_j$ and
 $x$.
Let $\tilde{h}$ denote the first line of $M_j$
 such that $x\in \tilde{h}^-$, and pick any line $h\in\AA$ that contains $x$ and is not
 parallel to $\tilde{h}$.  Let $y:=h\cap\tilde{h}$.
By construction
 $h\in\bigcup_{i>j}M_i$, and  
since
 $x$ is between $y$ and $h\cap \ell$ on $h$, by Remark
 \ref{inriga} there is no ordering that is special w.r.t. $\ell$ and
 in which $y$ comes after $x$, as recursive orderability with respect
 to $\ell$ would
 require.   


On the other hand, if $\AA$ {\em is} complete at every $a_j$, then an
explicit recursive combinatorial polar ordering can be described as
follows. Write $\AA=\{h_0,h_1,\ldots\}$ according to the order in which
the lines intersect $\ell$. To begin with, being complete implies that there every point
contained in $h_0^-$ lies actually on $h_0$. It is now evident that
the sequence $a_1,a_2,\ldots$ is a valid sequence of flippings, that
leads to a pseudoline $\ell_1$ with every point in $P\cap h_0$ on its
``backside''. Because there are no points in the interior of the cone
$h_1^+\cap h_2^-$, clearly one can now perform the flips across all
points of $h_2$. Clearly one can go on this way until the moving
pseudoline has flipped across every point in $\Lambda_1$.

We leave it to the reader to check that now one can perform all the flips of
points in $\Lambda_j$ for increasing $j$, each time following the
order of lines induced by the intersection with $\ell$.\end{pfsk}

We obtain a complete characterization of recursively orderable arrangements in the plane.

\begin{thm}\label{thm_2d} An affine arrangement of lines in the plane is recursively orderable
  if and only if there is a general flag $(b,\ell)$ so that $\AA$ is complete with respect to $(b,\ell)$.
\end{thm}

Some general facts about recursively orderable arrangements can be deduced.
\begin{rem} {\em Not all real reflection arrangements are recursively orderable.} For
  example consider the arrangement of type $H_3$. This is a central
  arrangement in $\R^3$, so it is recursively orderable if and only if there is a
  generic section of it that is recursively orderable. If we consider the
  projection of the associated dodecahedron on the plane of the
  section, we see that the points of this arrangement of lines
  correspond to vertices, to centers of edges or to centers of
  pentagonal faces. It is easy to see by case-by-case inspection that for
  every choice of $a_0$, of an adjacent chamber as $B$ and of a
  suitable line for $\ell$, $\Lambda_1$
  is never complete with respect to $(b,\ell)$. Indeed, if $a_0$
  corresponds to a pentagon $p$, the obstruction comes from a point
  corresponding to an edge $e$ that is not adjacent to $p$ but belongs to a pentagon adjacent to $p$ (and vice-versa), while the obstruction for
  every `vertex-type' choice of $a_0$ comes from another vertex that
  belongs to a common pentagon, but is not adjacent to $a_0$.
\end{rem}

\begin{rem} {\em Not all recursively orderable arrangements are $K(\pi,1)$.} A
  counterexample can in fact be given already in dimension $3$:
  consider the generic arrangement with defining form $xyz(x+y+z)$ in $\mathbb{R}^3$. By Hattori's theorem, this arrangement is not aspherical (see \cite[Corollary 5.23]{OT}). However, it is central and any $2$-dimensional section of it is easily seen to be recursively orderable.
\end{rem}

\section{Supersolvable arrangement are recursively orderable.}\label{supersolvable}

The class of ``strictly linearly fibered'' arrangements was introduced by Falk and Randell \cite{FaRa} in order to generalize the technique of Fadell and Neuwirth's proof \cite{FaNe} of asphericity of the braid arrangement (involving a chain of fibrations). Later on, Terao \cite{terao1} recognized that strictly linearly fibered arrangements are exactly those which intersection lattice is { supersolvable} \cite{stanley1}. Since then these are known as {\em supersolvable arrangements}, and deserved intense consideration. 

The goal of this section is to prove that every supersolvable real arrangement is recursively orderable. Let us begin by the definition.

\newcommand{\rank}{\textup{rank}}
\begin{df}\label{supers} A central arrangement $\AA$ of complex hyperplanes in $\mathbb{C}^d$ is called {supersolvable}
if there is a filtration
$\AA=\AA_d \supset \AA_{d-1} \supset \cdots \supset \AA_2 \supset \AA_1$
such that
\begin{enumerate}
\item[(1)] $\rank(\AA_i)=i$ for all $i=1,\ldots ,d$
\item[(2)] for every two $H, H^{\prime} \in \AA_i$ there exits
some $H^{\prime \prime} \in \AA_{i-1}$ such that 
$H \cap H^{\prime} \subset H^{\prime \prime}$.
\end{enumerate}
\end{df}

Before getting to the actual theorem, let us point out the key geometric fact.

\begin{rem}\label{generalext} Let $\AA$ be as in Definition \ref{supers} and consider the arrangement $\AA_{d-1}$ in $\mathbb{R}^d$. It is clearly not essential, and the top element of $\LL(\AA_{d-1})$ is a 1-dimensional line that we may suppose to coincide with the $x_1$-axis. The arrangement $\AA_{d-1}$ determines an essential arrangement on any hyperplane $H$ that meets the $x_1$-axis at some $x_1=t$. 
For all $t$, the intersection of $\AA_{d-1}$ with the hyperplane $H$ determines an essential, supersolvable arrangement $\AA'_{d-1} \subset \mathbb{R}^d$ with $\AA'_{r}=\AA_r$ as sets, for all $r\leq d-1$. Thus, given a flag of general position subspaces for $\AA'_{d-1}$, we can find a combinatorially equivalent flag $(V_k)_{k=0,\ldots, d-2}$ on  $H$.

Now let us consider a hyperplane $H$ in $\mathbb{R}^d$ that is orthogonal to the $x_1$-axis, and suppose we are given on it as above a valid flag $(V_k)_{k=0,\ldots,d-2}$ of general position subspaces for $\AA_{d-1}$. By tilting $H$ around $V_{d-2}$ 
we can  obtain a hyperplane $H'$ that is in general position with respect to $\AA$ and for which all points of $\AA\cap H'$ are on the same side with respect to $V_{d-2}$, and for which $V_0$ lies in an unbounded chamber.

By setting $V_{d-1}:=H'$, $V_d:=\mathbb{R}^d$ we thus obtain a valid general flag for $\AA=\AA_d$. Define $\PP^k(\AA_d)$ as the points of $\AA_d\cap V_k$ and analogously for $\PP^k(\AA_{d-1})$. The flag remains general by translating $H'=V_{d-1}$ in $x_1$-direction away from the origin: we can therefore suppose that there is $R\in\mathbb{R}$ such that  for all $k$, $k=1,\ldots,d-1$, every element of $\PP^k(\AA_{d-1})$ is contained in a ball of radius $R$ centered in $V_0$, that contains no element of $\PP^k(\AA_{d})\setminus\PP^{k}(\AA_{d-1})$.

\end{rem}

\begin{crl}\label{segmentato}
Let $\AA$ and $(V_k)_{k=1,\ldots,d}$ be as in the construction of Remark \ref{generalext}. Then, for every $k=1,\ldots,d$, if $F_1\in\PP^k(\AA_{d-1})$ and $F_2\in\PP^k(\AA)\setminus \PP^{k}(\AA_{d-1})$ are both contained in the support of the same $F\in\PP^{k-1}(\AA)$, then $F_1\leadsto^k F_2$ in every special ordering of $\PP^k(\AA)$.
\end{crl}
\begin{pf} This is an immediate consequence of Remark \ref{inriga} and \ref{generalext}.
\qed\end{pf}

\begin{thm}\label{ssfol} Any supersolvable complexified arrangement $\AA$ is recursively orderable. 
Moreover, the recursively orderable special ordering $\leadsto$ can be chosen so that
for all $i=2,\ldots,d$ and all $k=1,\ldots, i-1$, if  $F_1\in\PP^k(\AA_{i-1})$ and  $F_2\in\PP^k(\AA_{i}) \setminus \PP^k(\AA_{i-1})$  lie in the support of the same $k+1$-codimensional face,
then $F_1 {\leadsto} F_2$.
\end{thm}

\begin{pf} If $\AA$ has rank one, there is nothing to prove. So let $d:=\rank(\AA)>1$ and suppose the claim holds for all complexified supersolvable arrangements or rank strictly less than $d$ - in particular, for $\AA_{d-1}$. 

The general flag $(V_k)_{k=0,\ldots,d}$ we will use is obtained via Remark \ref{generalext} from a general flag for $\AA_{d-1}$ that gives rise to a special ordering satisfying the claim of the theorem. In particular, there exists a special ordering of $\PP(\AA_{d-1})$ that satisfies the property required by the claim for every $i=2,\ldots,d-2$ (and every $k=0,\ldots,i-1$). By Corollary \ref{segmentato}
and Remark \ref{trovarli}, we only have to describe, for every $k$, a special ordering  of $\PP^k(\AA)$ that satisfies condition (ii)(a) of Definition \ref{df:maindef}. This will be done by a new induction on $k$.

For $k=0$ there is nothing to prove, and for $k=1$ the only possible special ordering will clearly do. Let then $k>1$. Suppose that  recursive special orderings $\leadsto^{k-2}, \leadsto^{k-1}$ have already been defined on $\PP^{k-2}$ and $\PP^{k-1}$, and write $\PP^{k-1}=\{p_1,p_2,\ldots\}$ accordingly. Since $\AA$ is supersolvable, every $F\in\PP^k(\AA)$ is contained in the support of some element of $\PP^{k-1}(\AA_{d-1})$ that we will call $p(F)$. So what we have to show is the following.
\begin{claim} The ordering on $\PP^k(\AA)$ defined by 
$$F_1\leadsto F_2\Leftrightarrow \left\{\begin{array}{l}p(F_1)\leadsto^{k-1}p(F_2)\textrm{ or }\\ p(F_1)=p(F_2)\textrm{ and }F_1\textrm{ is between }p(F_2)\textrm{ and }F_2\textrm{ on }\vert p( F_2) \vert \end{array}\right.$$ is a special ordering.
\end{claim}
\begin{pfcl} Consider a special ordering of $\PP^k(\AA)$ that agrees with the above ordering up to some face $F_1$, and suppose for contradiction that $F_1$ is not near the moving pseudohyperplane, i.e., that there is $F_2$ with $p(F_1)\leadsto^{k-1} p(F_2)$ which is on a line passing through $F_1$ between $F_1$ and the moving pseudohyperplane. By the inductive hypothesis on $\AA_{d-1}$ we know that the above defined ordering is indeed special for the elements of $\PP^k(\AA_{d-1})$, and by Corollary \ref{segmentato} we conclude that $F_1$ cannot be in $\PP(\AA_{d-1})$.

Thus, the only obstruction to the construction of such a total ordering would come from the following situation: two faces $F_1,F_2\in\PP^k(\AA)\setminus\PP^k(\AA_{d-1})$ lying on the support of the same $q\in\PP^{k-1}(\AA)\setminus \PP^{k-1}(\AA_{d-1})$ so that $p(F_1)\leadsto^{k-1}p(F_2)$ but $F_2$ lies between $q$ and $F_1$ on $\vert q \vert$. We prove that this situation can indeed not occur.

\begin{figure}[h]\begin{center}\begin{tabular}{cc}
\includegraphics[scale=0.4]{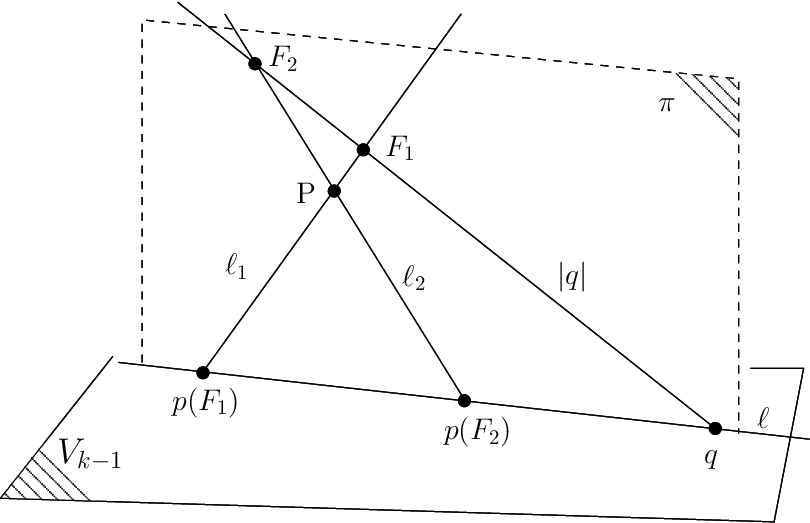}&
\includegraphics[scale=0.4]{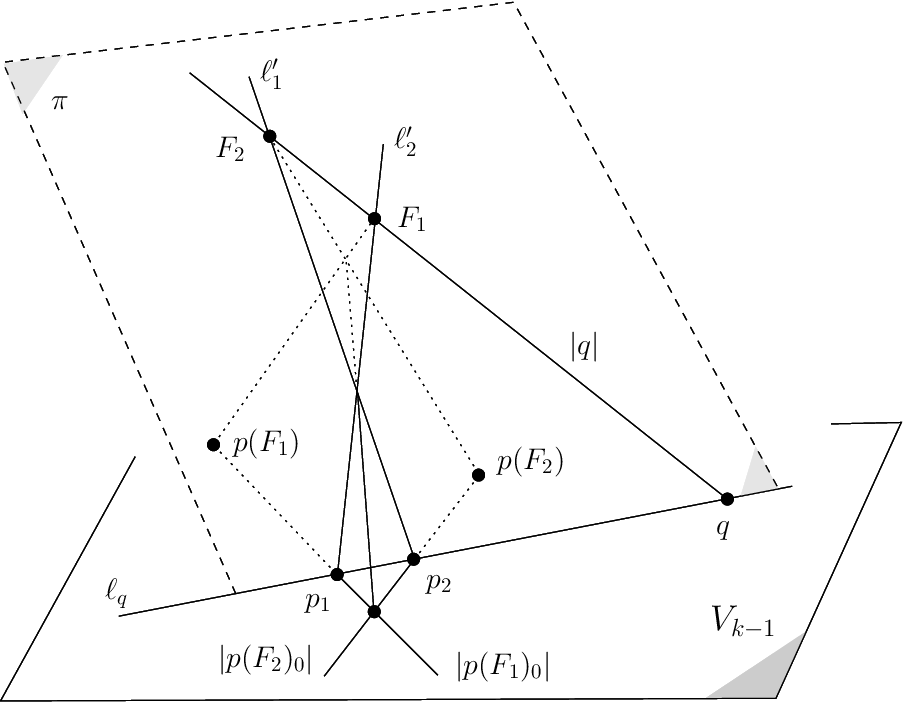}\\
(1) & (2)
\end{tabular}\end{center}
\caption{$\,$}\label{duefigure}
\end{figure}

 Given any $p\in\PP^{k-1}(\AA)$, let $p_0:=\min\{x\in\PP^{k-2}(\AA)\mid p\subset \vert x \vert\}$ as in Definition \ref{cinque}. Then we have two cases.

 {\em Case 1} (see Figure \ref{duefigure}.(1)) $p(F_1)_0=p(F_2)_0$. This means $p(F_1),p(F_2)\in\ell$, where $\ell:=\vert p(F_1)_0\vert$. The line $\ell$ is the intersection $\pi\cap V_{k-1}$ of $V_{k-1}$ with a plane $\pi$ in $V_k$ that contains also the lines $\ell_1:=\vert p(F_1)\vert$ and $\ell_2:=\vert p(F_2)\vert$. Then this plane must contain also the line $\vert q \vert$. Since $\AA_{d-1}$ is central, $\ell_1$ and $\ell_2$ must intersect, and this gives a point $P\in\PP^k(\AA_{d-1})$ that, by Remark \ref{inriga}, lies between $p(F_i)$ and $F_i$ for $i=1,2$. Again, by Remark \ref{inriga} we know that on $\ell$ we have the sequence of points $q,p(F_2),p(F_1)$, so on $|q|$ we have the sequence $q, F_1 , F_2$, and there is no obstruction.

 {\em Case 2} (see Figure \ref{duefigure}.(2)). $p(F_1)_0\leadsto p(F_2)_0$. Since $q\in \PP(\AA)\setminus \PP(\AA_{d-1})$, as above we have that the line $\ell_q:=\vert q_0 \vert$ intersects $\vert p(F_i)_0\vert$ in a point $p_i$ between $p(F_i)$ and $p(F_i)_0$, for $i=1,2$. Consider now the plane $\pi$ spanned by $\vert q \vert$ and $\ell_q$ (this might not be a flat of $\AA$), and on it, for $i=1,2$ the line $\ell'_i$ spanned by $p_i$ and $F_i$. The intersection $\ell'_1\cap\ell'_2$ lies on the segments $\overline{p_1F_1}$ and $\overline{p_2F_2}$  only if $\vert p(F_1)_0\vert \cap \vert p(F_2)_0 \vert$ is between $p(F_i)_0$ and $p_i$ Since the Theorem holds in $V_{k-1}$ it is now a straightforward check to verify that $p(F_1)\leadsto p(F_2)$ implies that $F_1$ lies between $F_2$ and $q$ on $\vert q \vert$ (Figure \ref{duefigure}.(2) describes one of the two possible cases - namely, when $\overline{p_1F_1}\cap \overline{p_2F_2}$ is not empty).  
\end{pfcl}
This concludes the proof of Theorem \ref{ssfol}. \qed\end{pf}

\bibliographystyle{elsarticle-num}
\bibliography{bibfollow}

\begin{thebibliography}{10}
\expandafter\ifx\csname url\endcsname\relax
  \def\url#1{\texttt{#1}}\fi
\expandafter\ifx\csname urlprefix\endcsname\relax\def\urlprefix{URL }\fi
\expandafter\ifx\csname href\endcsname\relax
  \def\href#1#2{#2} \def\path#1{#1}\fi


\bibitem{BLSWZ}
A.~Bj{\"o}rner, M.~Las~Vergnas, B.~Sturmfels, N.~White, G.~M. Ziegler, Oriented
  matroids, 2nd Edition, Vol.~46 of Encyclopedia of Mathematics and its
  Applications, Cambridge University Press, Cambridge, 1999.


\bibitem{delu}
E.~Delucchi, Shelling-type orderings of regular {CW}-complexes and acyclic
  matchings of the {S}alvetti complex, Int. Math. Res. Not. IMRN (2008) Art. ID
  rnm167, 39.




\bibitem{dimcapapa}
A.~Dimca, S.~Papadima, Hypersurface complements, {M}ilnor fibers and higher
  homotopy groups of arrangments, Ann. of Math. (2) 158 (2003) 473--507.


\bibitem{Edels}
H.~Edelsbrunner, L.~J. Guibas, Topologically sweeping an arrangement, J.
  Comput. System Sci. 38 (1989) 165--194, 18th Annual ACM Symposium on Theory
  of Computing (Berkeley, CA, 1986).

\bibitem{FaNe}
E.~Fadell, L.~Neuwirth, Configuration spaces, Math. Scand. 10 (1962) 111--118.

\bibitem{FaRa}
M.~Falk, R.~Randell, The lower central series of a fiber-type arrangement,
  Invent. Math. 82 (1985) 77--88.

\bibitem{FoLa}
J.~Folkman, J.~Lawrence, Oriented matroids, J. Combin. Theory Ser. B 25 (1978)
  199--236.

\bibitem{ForDM}
R.~Forman, Morse theory for cell complexes, Adv. Math. 134 (1998) 90--145.


\bibitem{FuTa}
K.~Fukuda, A.~Tamura, Local deformation and orientation transformation in
  oriented matroids, Ars Combin. 25 (1988) 243--258, eleventh British
  Combinatorial Conference (London, 1987).


\bibitem{Kozlov}
D.~Kozlov, Combinatorial algebraic topology, Springer, Berlin, 2008.


\bibitem{OT}
P.~Orlik, H.~Terao, Arrangements of hyperplanes, Springer, Berlin, 1992.

\bibitem{randell}
R.~Randell, Morse theory, {M}ilnor fibers and minimality of hyperplane
  arrangements, Proc. Amer. Math. Soc. 130 (2002) 2737--2743 (electronic).


\bibitem{JRG}
J.~Richter-Gebert, Testing orientability for matroids is {NP}-complete, Adv. in
  Appl. Math. 23~(1) (1999) 78--90.


\bibitem{salvetti}
M.~Salvetti, Topology of the complement of real hyperplanes in {${\bf C}\sp
  N$}, Invent. Math. 88 (1987) 603--618.

\bibitem{salvettisette}
M.~Salvetti, S.~Settepanella, Combinatorial {M}orse theory and minimality of
  hyperplane arrangements, Geom. Topol. 11 (2007) 1733--1766.


\bibitem{stanley1}
R.~P. Stanley, Supersolvable lattices, Algebra Universalis 2 (1972) 197--217.

\bibitem{terao1}
H.~Terao, Modular elements of lattices and topological fibration, Adv. in Math.
  62 (1986) 135--154.

\end{thebibliography}

\end{document}